\input ./style/arxiv-general.cfg
\documentclass[MSNbibl,number,citesort,dvips]{arxbj}
\makeatletter
   \@ifpackageloaded{graphicx}{}{\usepackage{graphicx}}
\makeatother
\usepackage{upgreek}

\volume{22}
\issue{2}
\pubyear{2016}
\firstpage{1026}
\lastpage{1054}
\doi{10.3150/14-BEJ685}
\docsubty{FLA}

\makeatletter
\newcommand{\rrvert}{\vert}
\newcommand{\llvert}{\vert}
\renewcommand{\mid}{|}
\newcommand{\Z}{\mathbb{Z}}
\newcommand{\N}{\mathbb{N}}
\newcommand{\R}{\mathbb{R}}
\renewcommand{\i}{\mathfrak{i}}
\renewcommand{\d}{\mathrm{d}}
\def\Var{\operatorname{Var}}
\def\EXP{\mathbb{E}}
\def\PROB{\mathbb{P}}
\def\E{\mathbb{E}}
\def\p{\mathbb{P}}
\newtheorem{lemma}{Lemma}
\newtheorem{theorem}{Theorem}
\newtheorem{proposition}{Proposition}
\newtheorem{corollary}{Corollary}
\newremark{remark}{Remark}
\makeatother

\begin{document}
\begin{frontmatter}

\title{Asymptotic behavior of the generalized St.~Petersburg sum
conditioned on its~maximum}
\runtitle{St.~Petersburg sum conditioned on its maximum}

\begin{aug}
\author[A]{\inits{G.}\fnms{G\'abor}~\snm{Fukker}\thanksref{A,e1}\ead[label=e1,mark]{fukkerg@math.bme.hu}},
\author[A]{\inits{L.}\fnms{L\'aszl\'o}~\snm{Gy\"orfi}\thanksref{A,e2}\ead[label=e2,mark]{gyorfi@cs.bme.hu}}
\and
\author[B]{\inits{P.}\fnms{P\'eter}~\snm{Kevei}\corref{}\thanksref{B}\ead[label=e3]{kevei@math.u-szeged.hu}}
\dedicated{Dedicated to the memory of S\'andor Cs\"org\H{o}}
\address[A]{Department of Computer Science and Information Theory,
Budapest University of Technology and Economics, 1521 Stoczek u. 2,
Budapest, Hungary.\\ \printead{e1,e2}}
\address[B]{MTA--SZTE Analysis and Stochastics Research Group, Bolyai
Institute, University of Szeged, 6720 Aradi
v\'ertan\'uk tere 1, Szeged, Hungary.
\printead{e3}}
\end{aug}

%
\received{\smonth{8} \syear{2013}}
%
\revised{\smonth{9} \syear{2014}}

%
\begin{abstract}
In this paper, we revisit the classical results on the
generalized St.~Petersburg sums. We determine the limit distribution of
the St.~Petersburg sum conditioning on its maximum, and we analyze how
the limit depends on the value of the maximum.
As an application, we obtain an infinite sum representation of the
distribution function of the possible semistable limits. In the representation,
each term corresponds to a given maximum, in particular this result
explains that the
semistable behavior is caused by the typical values of the maximum.
\end{abstract}

%
\end{frontmatter}

\section{Introduction} \label{sectintro}

Peter offers to let Paul toss a possibly biased coin repeatedly
until it lands heads and pays him $r^{k/\alpha}$ ducats if this
happens on the $k$th toss, where $k\in\N= \{1,2,\ldots\}$,
$p \in(0,1)$ is the probability of heads at each throw, $q=1-p$,
$r=q^{-1}$, while $\alpha> 0$ is a payoff parameter.
This is the so-called generalized St.~Petersburg
game with parameter $(\alpha, p)$. The classical St.~Petersburg game
corresponds to $\alpha= 1$ and $p = 1/2$.
If $X$ denotes Paul's winning in this St.~$\operatorname{Petersburg}(\alpha, p)$ game, then
$\PROB\{ X = r^{k/\alpha} \} = q^{k-1} p$, $k\in\N$.
Put $\lfloor x \rfloor$ for the lower integer part,
$\lceil x \rceil$ for the upper integer part and
$\{ x \}$ for the fractional part of $x$. Then
the distribution function of the gain is
%
%
\begin{equation}
\label{eqdist-stp} F (x) = \PROB\{ X \leq x \} = \cases{ 0, &\quad$x <
r^{1/\alpha}$,
\hspace*{1pt}\cr
\displaystyle1- q^{\lfloor\alpha \log_r x \rfloor} = 1- \frac{ r^{
\{ \alpha\log_r x \} } }{x^\alpha},
&\quad$x \geq r^{1/\alpha}$,}
\end{equation}
where $\log_r$ stands for the logarithm to the base $r$.

In the following all the functions, constants and random variables depend
on the parameters $\alpha$ and $p$. For the sake of readability
we suppress everywhere the upper index ${\alpha, p}$.

We see that the
payoff parameter $\alpha>0$ is in fact a tail parameter of the distribution.
In particular, $\E(X^\alpha) = \infty$, but
$\E(X^\beta) = p/(q^{\beta/\alpha} - q)$
is finite for $\beta\in(0,\alpha)$, so for $\alpha> 2$ Paul's gain $X$ has
a finite variance, so L\'evy's central limit theorem holds.
As Cs\"org\H{o} pointed out in \cite{Cs02}, even\vadjust{\goodbreak} for $\alpha=2$ the
St.~$\operatorname{Petersburg}(2,p)$ distribution is in the domain of attraction of
the normal law. This can be checked by straightforward calculation,
using the well-known characterization of the~domain of attraction of
the normal
law.
Hence, the case $\alpha\ge2$ is substantially
different from the more difficult case $\alpha< 2$.
In Section~\ref{sectcond},
when we are dealing with asymptotic behavior of the sums as $n \to
\infty$
we usually assume that $\alpha< 2$. We indicate the possible values
of $\alpha$ in all of the~statements. Of course, the most interesting
case is the classical one, when
$\alpha= 1$, for which the mean is infinite.

\subsection{The sum}

Let $X, X_1, X_2, \ldots$ be i.i.d. St.~$\operatorname{Petersburg}(\alpha, p)$ random
variables,
let $S_n = X_1 + \cdots+ X_n$ denote their partial sum, and
$X_n^* = \max_{1 \leq i \leq n} X_i$ their maximum.
Since the bounded oscillating function $r^{\{ \alpha \log_r x \} }$ in
the numerator
of the distribution function in
(\ref{eqdist-stp}) is not slowly varying at infinity, by the classical
Doeblin--Gnedenko criterion (cf. \cite{GK}) the underlying St.~Petersburg
distribution is \textit{not} in the domain of attraction of any stable
law. That is there is
no asymptotic distribution for $(S_n - c_n)/a_n$, in the usual sense,
whatever the
centering and norming constants are. This is where the main difficulty lies
in analyzing the St.~Petersburg games.

However, asymptotic distributions do exist along subsequences
of the natural numbers. In the classical case, when $\alpha= 1$, $p=1/2$,
Martin-L\"of \cite{ML} ``clarified the St.~Petersburg
paradox,'' showing that $S_{2^k}/2^k - k$ converges in distribution, as
$k \to\infty$.
Cs\"org\H{o} and Dodunekova \cite{CsD} showed that there are continuum
of different
types of asymptotic
distributions of $S_n/n - \log_2 n$ along different subsequences of
$\mathbb{N}$.

In order to state the necessary and sufficient condition for the
existence of the limit,
we introduce the positional parameter
%
%
\begin{equation}
\label{eqgamma-def} \gamma_n = \frac{n}{r^{\lceil\log_r n \rceil}}
\in(q, 1],
\end{equation}
which shows the position of $n$ between two consecutive powers of $r$.
Put
%
%
\begin{equation}
\label{eqmu} \mu_n = \cases{ \displaystyle n^{1- \alpha^{-1}} \frac{p}{q^{1/\alpha}
- q}, &
\quad for $\alpha\ne1$,
\vspace*{5pt}\cr
\displaystyle\frac{p}{q} \log_r
n, &\quad for $\alpha= 1$.}
\end{equation}
In Theorem 1 in \cite{Cs02}, Cs\"org\H{o} showed that the following
merging theorem holds
(in fact a sharp estimate for the rate is also provided):
%
%
\begin{equation}
\label{eqsum-merge} \sup_{x \in\mathbb{R}} \biggl\llvert\PROB\biggl
\{
\frac{S_n}{n^{1/\alpha}} - \mu_n \leq x \biggr\} - G_{\gamma_n} (x)
\biggr\rrvert\to0\qquad\mbox{as } n \to\infty,
\end{equation}
where $G_\gamma$ is the distribution function of the infinitely
divisible random
variable $W_\gamma$, $\gamma\in(q, 1]$ with
characteristic function
%
%
\begin{equation}
\label{eqy} \EXP\bigl( \mathrm{e}^{ \i t W_{\gamma} } \bigr) 
= \mathrm{e}^{y_\gamma(t) }
= \exp\biggl( \i t [ s_\gamma+ u_{\gamma} ] + \int
_0^{\infty} \biggl( \mathrm{e}^{\i t x} - 1 -
\frac{\i t x }{1+x^2} \biggr) \,\d R_{\gamma} (x) \biggr)
\end{equation}
with
\begin{eqnarray*}
s_\gamma& =& \cases{\displaystyle\frac{p}{q - q^{1/\alpha}} \frac{1}{\gamma^{(1-\alpha
)/\alpha}},&\quad  $\alpha\ne1$,
\vspace*{5pt}\cr
\displaystyle \frac{p}{q} \log_r \frac{1}{\gamma}, &\quad $\alpha= 1$,}
\\
u_\gamma& =& \frac{p}{q} \gamma^{(\alpha+1)/\alpha} \sum
_{k=1}^\infty\frac{r^{((1-\alpha)/\alpha) k}}{\gamma^{2/\alpha} + r^{2k/\alpha}} - \frac{p}{q}
\gamma^{(\alpha-1)/\alpha}\sum_{k=0}^\infty
\frac{1}{\gamma^{2/\alpha} r^{((3 - \alpha)/\alpha) k} +
r^{((1-\alpha)/\alpha) k}},
\end{eqnarray*}
and L\'evy function
%
%
\begin{equation}
\label{eqLevy-func} R_{\gamma} (x) = - \gamma q^{\lfloor\log_r
(\gamma x^\alpha) \rfloor} = -
\frac{r^{\{ \log_r (\gamma x^\alpha) \} } }{x^\alpha}, \qquad x>0.
\end{equation}
From this form, it is clear that $W_\gamma$ is a semistable random variable
with characteristic exponent~$\alpha$.
For the precise rate of the convergence in (\ref{eqsum-merge}) see Cs\"
org\H{o}
\cite{Cs07}, where short merging asymptotic expansions are provided, and
also additional historical background and references are given.
Merging asymptotic expansions are proved by Pap \cite{pap}, where the
length of the
expansion depends on the parameter $\alpha$: the closer $\alpha$ is to
0, the longer expansion
is possible. Pap \cite{pap} also shows non-uniform asymptotic expansions.
The natural framework of the merging theorems is the class of
semistable distributions,
see Cs\"org\H{o} and Megyesi \cite{CsMe02}. In Section~\ref{subsecttyp-max}, we
briefly collect the definition and basic properties of semistable distributions.

\subsection{The maximum}

It turns out that the maximum $X_n^*$ has similar asymptotic behavior
as the sum.
Let us consider the classical case again, that is, $\alpha= 1, p=1/2$.
For $\gamma\in(1/2,1]$, introduce the distribution function
\[
H_\gamma(x) = \cases{ 0, &\quad for $x \leq0$,
\cr
\displaystyle\exp
\bigl( - \gamma2^{-\lfloor\log_2 (\gamma x) \rfloor} \bigr), &\quad
for $x > 0$.}
\]
Berkes, Cs\'aki and Cs\"org\H{o} \cite{BCsCs} showed that although
there is no limit
theorem for the normed maximum through the whole sequence, the
following merging theorem holds:
%
%
\begin{equation}
\label{eqmax-merge} \sup_{x \in\R} \biggl\llvert\PROB\biggl\{
\frac{X_n^*}{n} \leq x \biggr\} - H_{\gamma_n}(x) \biggr\rrvert= \mathrm{O}\bigl(n^{-1}\bigr) \qquad\mbox{as } n \to\infty,
\end{equation}
with the positional parameter $\gamma_n$ defined in (\ref{eqgamma-def}).
Note that even though the ``limiting'' distribution function is not
continuous, merging holds in
uniform distance.
A more general setup is treated by Megyesi \cite{Megyesi}, see in particular
Theorem 4 in \cite{Megyesi}.

%
\begin{figure}[b]

\includegraphics{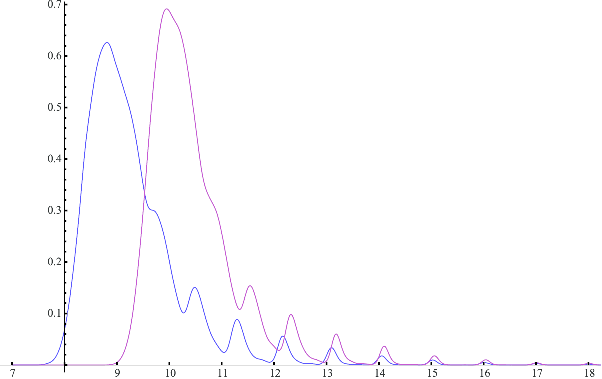}

\caption{The histograms of $\log_2 S_{n}$ for $n=2^{6}$ and for
$n=2^{7}$.}\vspace*{-8pt} \label{fig1}
\end{figure}

%
\begin{figure}

\includegraphics{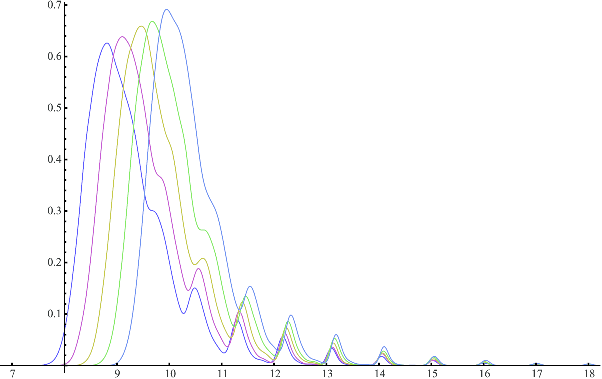}

\caption{The histograms of $\log_2 S_{n}$ for $n=2^{6+\eta}$,
$\eta=0, 0.25, 0.5, 0.75, 1$.}\vspace*{-8pt} \label{fig2}
\end{figure}

The merging theorems (\ref{eqsum-merge}) and (\ref{eqmax-merge})
immediately imply that in the classical case
$S_n/n - \log_2 n$ and $X_n^*/n$ converges along the subsequence
$\{ n_k \}$ if and only
if $\gamma_{n_k} \to\gamma$, as $k \to\infty$, for some $\gamma\in
[1/2, 1]$,
or $\{ \gamma_{n_k} \}$ has exactly
two limit points, $1/2$ and $1$. The latter is called
\emph{circular convergence}, as it can be seen as convergence
on the interval $[1/2, 1]$, $1/2$ and $1$ identified. See \cite{Cs02}
and \cite{Cs07}.
Similar statement holds in the general case.

Having seen these similarities it
is tempting to investigate the maximum and the sum together.
In Figures~\ref{fig1} and \ref{fig2} (all the figures correspond to the
classical
case), one can see that the histograms of $\log_2 S_n$ are mixtures of
unimodal densities
such that the first lobe is a mixture of overlapping densities,
while the side lobes have disjoint support.
For doubling $n$, in Figure~\ref{fig1} the pairs of
corresponding side lobes are almost identical, which suggests an oscillating
behavior governed by the parameter $\gamma_n$ in (\ref{eqgamma-def}).
Figure~\ref{fig2} shows the histograms of $\log_2 S_{n}$ for
$n=\lfloor2^{6+\eta} \rfloor$, $\eta=0, 0.25, 0.5, 0.75, 1$, that is
for different
values of $\gamma_n$.

We mention that investigating the joint behavior of the sum and the maximum
goes back to Chow and Teugels \cite{chow-teugels}.
Let $Y, Y_1, Y_2, \ldots$ be i.i.d. random variables, $Z_n$ and $Y_n^*$
their partial sum and partial maximum, respectively.
In \cite{chow-teugels}, Chow and Teugels show that for some
deterministic sequences
$a_n > 0, c_n >0, b_n, d_n$,
$(Z_n/a_n - b_n, Y_n^* / c_n - d_n)$ converges in distribution to $(U,V)$,
where neither $U$ nor $V$ is
degenerate, if and only if $Y$ belongs to the domain of attraction of a
stable law, and also belongs to the maximum
domain of attraction of some extreme value distribution. Moreover, they
also characterize when $U$ and $V$ are
independent. The key technique in their proof is the ``hybrid''
function: characteristic
function of the sum, and distribution function of the maximum.
The same results using point process methods were proved by Kasahara
\cite{kasahar} and
by Resnick \cite{Resnick}.
Arov and Bobrov \cite{arovbobrov} consider the maximum modulus term
instead of the maximum.
The joint convergence
is also studied in case of non-independent random variables, we only
mention a recent
paper by Silvestrov and Teugels \cite{SilvestrovTeugels}.
Without the proof, we mention that the method of Chow and Teugels can
be used to obtain subsequential joint limit theorems for the sum and
for the maximum
in our setup.

%
\begin{figure}

\includegraphics{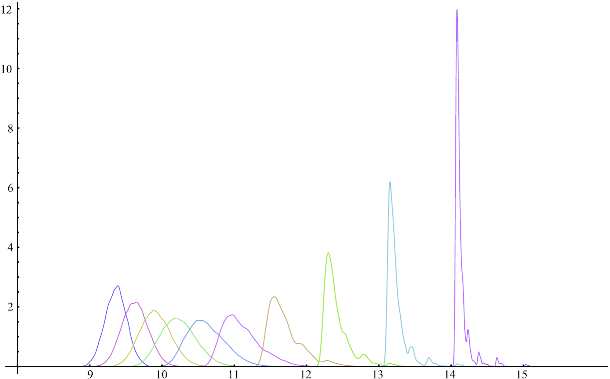}

\caption{The conditional histograms for $\log_2 S_{n}$, $n=2^7$.}\vspace*{-8pt}\label{fig3}
\end{figure}

In the present paper, we investigate together the maximum and the sum
of the St.~Petersburg random variables.
In Section~\ref{sectcond}, we determine the asymptotic distribution of $S_n$
conditioning on the maximum value, and we demonstrate how the limit
depends on the maximum.
Figure~\ref{fig3} shows the different blocks of the
smoothed histogram of $\log_2 S_n$, $n=2^7$, such that
in each block the maximum is the same, that is each lobe is the
smoothed conditional
histogram for $S_n$ given that $X_n^* = 2^{k}$, for $k=5,6, \ldots, 14$.
Comparing it with Figure~\ref{fig1} it is visible that the lobes are
determined by the
behavior of the maximum term. As (\ref{eqmax-merge}) states, the
typical value for $k$ is $\log_2 n$. The first
lobes correspond to smaller values of $X_n^*$, and so it is natural to
expect a Gaussian
limit; 
Proposition \ref{propnormality} deals with this case. The typical
values of the maximum
make the important contribution, and this is where the limiting
semistable law appears.
The middle lobes are the density functions of infinitely divisible
distribution functions,
each of these has finite expectation. This conditional limit theorem is stated
in Proposition \ref{proptypical-max-cond}.
Finally, as the maximum becomes larger and larger
it dominates the whole sum $S_n$. 
The conditional
limit for large maximum is contained in Proposition \ref{proplarge-max}.

In Section~\ref{sectappl}, we consider an application of this approach.
As a consequence of Proposition \ref{proptypical-max-cond}, in Theorem
\ref{thmixture}
we show that
\[
G_{\gamma} (x) = \sum_{j=-\infty}^\infty
\widetilde G_{j, \gamma} (x) p_{j,\gamma},
\]
where $G_\gamma$ is the merging distribution function appearing in (\ref
{eqsum-merge}).
Here $\widetilde G_{j,\gamma_n}$ corresponds to the distribution
function of the sum conditioned
on $X_n^* = r^{(\lceil\log_r n \rceil+ j)/\alpha}$, and $p_{j, \gamma
_n}$ is the approximate probability
of this event. The decomposition shows that the merging property is
caused by the
asymptotic properties of the maximum.

Finally, we note that recently Gut and Martin-L\"of \cite{GML}
investigated the
so-called max-trimmed St.~Petersburg games in the classical case,
where from the sum all the maximal observations are discarded.
They obtained the asymptotic behavior of the trimmed sum
along subsequences of the form
$(\lfloor\gamma2^n \rfloor)_{n \in\N}$.

\section{Conditioning on the maximum} \label{sectcond}

In this section, first we revisit the limit properties of $X_n^*$, and
then conditioning on different values of the maximum, we determine
the limit distribution of the sums.

\subsection{Asymptotics of the maximum}

For $j \in\Z$ and $\gamma\in[q, 1]$ introduce the notation
%
%
\[
\label{eqp} p _{j,\gamma} = \mathrm{e}^{- \gamma q^{ j}} \bigl( 1- \mathrm{e}^{- \gamma
(r-1) q^{j}}
\bigr).
\]
The following lemma is a reformulation of (\ref{eqmax-merge}) in the
general case.
We give the short proof for completeness. Recall the definition of
$\gamma_n$ in (\ref{eqgamma-def}).

%
\begin{lemma} \label{lemmamax-conv}
For any $\alpha> 0$ we have that
%
%
\begin{equation}
\label{eqmax-pr-merge} \sup_{j \in\Z} \bigl\llvert\PROB\bigl\{
X_n^* = r^{(\lceil\log_r n \rceil+ j)/\alpha} \bigr\} - p_{j,
\gamma_n} \bigr\rrvert= \mathrm{O}
\bigl(n^{-1}\bigr).
\end{equation}
In particular for any $j \in\Z$, as $n \to\infty$
\[
\PROB\bigl\{ X_n^* = r^{(\lceil\log_r n \rceil+ j)/\alpha}
\bigr\} \sim
\mathrm{e}^{- \gamma_n q^{ j}} \bigl( 1- \mathrm{e}^{- \gamma_n (r-1) q^{j}} \bigr).
\]
%
\end{lemma}

\begin{pf}
For any $k =1,2, \ldots$ we have
$\PROB\{ X_n^* \leq r^{k/\alpha} \} = ( 1 - q^{k} )^n$, and so
\begin{eqnarray*}
\bigl\llvert\PROB\bigl\{ X_n^* \leq r^{(\lceil\log_r n \rceil+ j)/\alpha} \bigr\} -
\mathrm{e}^{-\gamma_n q^{j}} \bigr\rrvert& =& \bigl\llvert\bigl( 1 - q^{ \lceil
\log_r n \rceil+ j }
\bigr)^n - \mathrm{e}^{-\gamma_n q^{j}} \bigr\rrvert
\\
& =& \biggl\llvert\biggl( 1 - \frac{\gamma_n q^{j}} {n} \biggr)^n -
\mathrm{e}^{-\gamma_n q^{j}} \biggr\rrvert
\\
& =& \mathrm{O}\bigl(n^{-1}\bigr).
\end{eqnarray*}
Since the latter holds uniformly, that is,
\[
\sup_{0 \leq y \leq n q} \biggl\llvert\biggl( 1 - \frac{y}{n}
\biggr)^n - \mathrm{e}^{-y} \biggr\rrvert= \mathrm{O}\bigl(n^{-1}
\bigr),
\]
and
\[
\PROB\bigl\{ X_n^* = r^{k/\alpha} \bigr\} = \PROB\bigl\{
X_n^* \leq r^{k/\alpha} \bigr\} - \PROB\bigl\{ X_n^*
\leq r^{(k-1)/\alpha} \bigr\},
\]
the proof is complete.
\end{pf}

%
\begin{remark}
The random variables $\alpha\log_r X_n^* -\lceil\log_r n \rceil$ have
a limit distribution
along subsequences $ \{ n_k = \lfloor\gamma r^{k} \rfloor\}_{k \in\N}$,
with $q < \gamma\leq1$, since using Lemma \ref{lemmamax-conv} above,
as $k \to\infty$
%
%
\begin{equation}
\label{eqlog-max-conv} \PROB\bigl\{ \alpha\log_r X_{n_k}^* -
\lceil\log_r n_k \rceil= j \bigr\} \to
\mathrm{e}^{- \gamma q^{ j}} \bigl( 1- \mathrm{e}^{- \gamma(r-1) q^{ j}} \bigr)
=p_{j,\gamma}.
\end{equation}

%
\begin{table}
\tabcolsep=0pt
\caption{Limit distribution of $\log_2 X_{n_{k}}^* -\lceil\log_2 n_{k}
\rceil$ in the classical case with $\gamma=1$}\label{tab1}
\begin{tabular*}{\tablewidth}{@{\extracolsep{\fill}}@{}lllllllll@{}}
\hline
$j$ & $-2$ & $-1$ & 0 & 1 & 2 & 3 & 4 & 5\\
\hline
$p_{j,1}$ & 0.018 & 0.117 & 0.233 & 0.239 & 0.172 & 0.104 & 0.057 & 0.03\\
\hline
\end{tabular*}
\end{table}
Table~\ref{tab1} contains the few largest values of $p_{j,1}$.
This is the main part of the limit distribution, as
$\sum_{j=-2}^{5} p_{j,1}\approx0.943$.

The asymptotic distribution (\ref{eqmax-pr-merge}) implies that
$\inf_n \Var(\log_r X_n^*)>0$, while in the classical case Gy\"orfi and Kevei
(Remark 2 in \cite{GyK}) showed that $\Var(\log_2 S_n)=\mathrm{O}(1/\log_2 n)$.
\end{remark}

%
\begin{remark}
Consider again the classical case.
We note that the merging theorem (\ref{eqlog-max-conv}) already
appears in F\"oldes \cite{foldes}.
Let $\mu(n)$ be the longest tail-run after tossing a fair coin $n$
times. Then Theorem 4 in
\cite{foldes} states that for any integer $j$
\[
\PROB\bigl\{ \mu(n) - \lfloor\log_2 n \rfloor< j \bigr\} =
\mathrm{e}^{-2^{-(j+1 - \{ \log_2 n \} )}} + \mathrm{o}(1).
\]

Since each single St.~Petersburg game lasts till to the first heads, in
our setup we are tossing the
coin until a random time, until heads appears $n$ times. Thus, the
number of tosses has a negative
binomial distribution with parameter $n$. Moreover, the values $(\log_2
X_k) -1$, $k=1,2, \ldots, n$,
are the number of tails between two consecutive heads, therefore the
quantity $\log_2 X_n^* - 1$ can be thought as
the longest tail-run in this coin tossing sequence.
\end{remark}

We investigate the conditional distribution of $S_n$ given that $X_n^*
= r^{k/\alpha}$.
The following lemma determines this conditional distribution. The
statement for
continuous random variables is much simpler, as in that case the
maximum value is
almost surely unique, and so $M_n = 1$ a.s. (see the definition below).
For the continuous version, see Lemma 2.1 in \cite{Darling}.

%
\begin{lemma} \label{lemmacond}
Let $Y, Y_1, \ldots, Y_n$ be discrete i.i.d. random variables with
possible values
$\{ y_1, y_2, \ldots\}$, $y_1 < y_2 < \cdots.$ Put
\[
G_k (y) = \p\{ Y \leq y \mid Y \leq y_k \}.
\]
Put $Z_n = Y_1 + \cdots+ Y_n$ for the partial sum, $Y_n^* = \max\{
Y_1, \ldots, Y_n \}$
for the partial maximum, and $M_n = | \{ k\dvt  1 \leq k \leq n, Y_k =
Y_n^* \} |$ for the
multiplicity of the maximum. Then given that
$Y_n^* = y_k$ and $M_n = m$
\[
Z_n \stackrel{\mathcal{D}} {=} m y_k +
Z_{n-m}^{(k-1)},
\]
where $Z_n^{(k-1)} = Y_1^{(k-1)} + \cdots+ Y_n^{(k-1)}$, with
$Y_1^{(k-1)}, \ldots, Y_n^{(k-1)}$
are i.i.d. with distribution function $G_{k-1}$.
\end{lemma}

\begin{pf}
We have
\begin{eqnarray*}
&&  \p\bigl\{ Z_n \leq y \mid Y_n^* = y_k,
M_n = m \bigr\}
\\
&&\quad = \frac{\p\{ Z_n \leq y, Y_n^* = y_k, M_n = m \}}{\p
\{ Y_n^* = y_k, M_n = m \}}
\\
&&\quad  = \frac{1}{\p\{ Y_n^* = y_k, M_n = m \}} \pmatrix{n
\cr
m}
\\
&&\qquad{} \times\p\Biggl\{ Y_1 = \cdots=
Y_m = y_k, \sum_{j={m+1}}^n
Y_j \leq y - m y_k, \max\{ Y_{m+1}, \ldots,
Y_n \} < y_{k} \Biggr\}
\\
&&\quad  = \frac{ {n\choose m} \p\{ Y = y_k \}^m \p\{ Y \leq y_{k-1} \}^{n-m}}{
\p\{ Y_n^* = y_k, M_n = m \}}
\\
&&\qquad{}  \times\p\Biggl\{ \sum_{j={m+1}}^n
Y_j \leq y - m y_k \mid\max\{ Y_j, j= m+1,
\ldots, n \} \leq y_{k-1} \Biggr\}
\\
&&\quad  = G_{k-1}^{*(n-m)} (y - m y_k),
\end{eqnarray*}
as stated.
\end{pf}

Put
%
%
\[
\label{eqN} N_n = \bigl\llvert\bigl\{ k\dvt  1 \leq k \leq n,
X_k = X_n^* \bigr\} \bigr\rrvert.
\]
According to the previous lemma in order to analyze the conditional
behavior of $S_n$, we first have to understand the
behavior of $N_n$.

%
\begin{lemma}
The conditional generating function of $N_n$ given $X_n^*$ is
%
%
\begin{equation}
\label{eqg-k,n} g_{k,n} (s) = \E\bigl[ s^{N_n} \mid
X_n^* = r^{k/\alpha} \bigr] = \frac{ ( 1 - q^{k-1} (1-ps) )^n - ( 1 -
q^{k-1} )^n }{
(1- q^k )^n - ( 1 - q^{k-1} )^n},
\end{equation}
and the generating function of $N_n$ is
%
%
\[
\label{eqg} g_n (s) = \E\bigl[ s^{N_n} \bigr] = \sum
_{k=1}^\infty\bigl[ \bigl( 1 -
q^{k-1} (1-ps) \bigr)^n - \bigl( 1 - q^{k-1}
\bigr)^n \bigr].
\]
\end{lemma}

\begin{pf}
Simply
%
\begin{eqnarray}
\label{eqNn-prob} \p\bigl\{ N_n = m \mid X_n^* =
r^{k/\alpha} \bigr\} & =& \frac{\p\{ N_n = m, X_n^* = r^{k/\alpha} \}
}{\p\{ X_n^* = r^{k/\alpha} \}}
\nonumber\\[-8pt]\\[-8pt]\nonumber
& =& \frac{{n\choose m} (q^{k-1} p)^m ( 1- q^{k-1})^{n-m} }{(1 - q^k)^n
- (1 - q^{k-1})^n}.
\end{eqnarray}
Therefore, by the binomial theorem the conditional generating function is
\begin{eqnarray*}
g_{k,n} (s) & =& \sum_{m=1}^n
s^m \frac{{n\choose m} (q^{k-1} p)^m ( 1- q^{k-1})^{n-m} }{(1 - q^k)^n
- (1 - q^{k-1})^n}
\\
& =& \frac{1}{(1 - q^k)^n - (1 - q^{k-1})^n} \bigl[ \bigl( s q^{k-1} p +
1 - q^{k-1}
\bigr)^n - \bigl(1 - q^{k-1}\bigr)^n \bigr].
\end{eqnarray*}
The unconditional version follows from the law of total probability.
\end{pf}

The distribution of $N_n$ in the classical case is calculated by
Gut and Martin-L\"of, in particular formula (\ref{eqNn-prob}) is
formula (4.1) in \cite{GML}. Moreover, in (4.3) in \cite{GML}
they determine the asymptotic behavior of $N_n$ conditioned on
typical maximum along geometric subsequences. This is formula~(\ref{eqNn-typ}) in the next proposition in the general
merging framework.

Now we can determine the asymptotic behavior of $N_n$.

%
\begin{proposition} \label{propNn}
Conditionally on $X_n^* = r^{k_n/\alpha}$, where $\log_r n -
k_n \to\infty$
%
%
\begin{equation}
\label{eqNn-norm} \frac{N_n - \E[ N_n \mid X_n^* = r^{k_n / \alpha}
] }{\sqrt{ \Var(N_n \mid X_n^* = r^{k_n/ \alpha} ) } } \stackrel
{\mathcal{D}} {\longrightarrow}
N (0,1)\qquad\mbox{as } n \to\infty.
\end{equation}
Conditionally on $X_n^* = r^{(\lceil\log_r n \rceil+ j)/{\alpha}}$, $j \in\Z$,
%
%
\begin{equation}
\label{eqNn-typ} \lim_{n \to\infty} \bigl\llvert g_{\lceil\log_r n
\rceil+ j, n} (s)
- h_{j,\gamma_n}(s) \bigr\rrvert= 0, \qquad s \in[0,1],
\end{equation}
where
%
%
\begin{equation}
\label{eqh} h_{j, \gamma} (s) = \frac{\mathrm{e}^{-(1-ps) \gamma q^{j-1} } -
\mathrm{e}^{-\gamma q^{j-1}} }{\mathrm{e}^{-\gamma q^j} - \mathrm{e}^{-\gamma q^{j-1}} },
\end{equation}
is the generating function of a Poisson$(pq^{j-1}\gamma)$ random
variable conditioned on not being zero.
While, if $k_n - \log_r n \to\infty$ then conditionally on $X_n^{*} =
r^{k_n / \alpha} $
%
%
\begin{equation}
\label{eqNn-large} N_n \stackrel{\p} {\longrightarrow} 1\qquad\mbox
{as }
n \to\infty.
\end{equation}
\end{proposition}

That is, we have three different regimes. In the typical range,
there are several random variables equal to the maximal value and the number
of these observations is distributed according to $h_{j,\gamma_n}$.
When the maximum is smaller than it should be, then there
are a lot of maximum values, while for too big values there is a single maximal
observation.

\begin{pf*}{Proof of Proposition \ref{propNn}}
Differentiating $g_{k,n}$ in (\ref{eqg-k,n}), we obtain
%
%
\begin{equation}
\label{eqNn-exp} \E\bigl[ N_n \mid X_n^* =
r^{k/\alpha}\bigr] = \frac{n q^{k-1} p (1 - q^k)^{-1}}{ 1 - (1 -
((p q^{k-1})/(1 - q^k)))^n }.
\end{equation}

First, we consider the case $\log_r n - k_n \to\infty$. Then
%
%
\begin{equation}
\label{eqaux1} \biggl(1 - \frac{p q^{k_n-1}}{1 - q^{k_n}} \biggr)^n
\to0,
\end{equation}
therefore
%
%
\begin{equation}
\label{eqNn-expasy} \E\bigl[ N_n \mid X_n^* =
r^{k_n/\alpha}\bigr] \sim\frac{n q^{k_n-1} p }{1 - q^{k_n}} =: \mu_{n,k_n}.
\end{equation}
(Note that we do not assume that $k_n \to\infty$ only that $\log_r n -
k_n \to\infty$.)
Using the simple identity that $\Var(N_n \mid X_n^* = r^{k/\alpha}) = g_{k,n}''(1) +
g_{k,n}'(1) - (g_{k,n}'(1))^2$, similar
computation gives
%
%
\begin{equation}
\label{eqNn-var} \Var\bigl(N_n \mid X_n^* =
r^{k_n/\alpha}\bigr) \sim\frac{n p q^{k_n-1}}{1 - q^{k_n}} \biggl( 1 -
\frac{p q^{k_n -1}}{1 - q^{k_n}}
\biggr) =: \sigma_{n,k_n}^2.
\end{equation}
Substituting into formula (\ref{eqg-k,n}), we have
\begin{eqnarray*}
&&\E\bigl[ \mathrm{e}^{\i t ({N_n - \mu_{n,k_n}})/{\sigma_{n,k_n}} } \mid X_n^*
= r^{k_n/ \alpha} \bigr]
\\
&&\quad  =
\mathrm{e}^{- \i t {\mu_{n,k_n}}/{\sigma_{n,k_n}} } \frac{ ( 1 - q^{k_n -1}
(1-p \mathrm{e}^{\i t/\sigma_{n,k_n} }) )^n - ( 1 - q^{k_n-1} )^n }{
(1- q^{k_n} )^n - ( 1 - q^{k_n-1} )^n}
\\
&&\quad  = \mathrm{e}^{- \i t {\mu_{n,k_n}}/{\sigma_{n,k_n}} } \frac{ ( 1 -
(p
q^{k_n -1} (1- \mathrm{e}^{\i t/\sigma_{n,k_n} })/(1- q^{k_n}) ))^n - ( 1 -
({p q^{k_n-1}}/({1 - q^{k_n}})) )^n }{
1 - ( 1 - ((p q^{k_n-1})/(1 - q^{k_n})) )^n }.
\end{eqnarray*}
By (\ref{eqaux1}), we have to determine the limit of
%
%
\begin{equation}
\label{eqBern} \mathrm{e}^{- \i t {\mu_{n,k_n}}/{\sigma_{n,k_n}} } \biggl(
1 - \frac{p q^{k_n -1} (1- \mathrm{e}^{\i t/\sigma_{n,k_n} })}{1- q^{k_n}}
\biggr)^n.
\end{equation}
Notice that
\[
1 - \frac{p q^{k_n -1} (1- \mathrm{e}^{\i t })}{1- q^{k_n}}
\]
is the characteristic function of a $0/1$ $\operatorname{Bernoulli}(pq^{k_n -1} / ( 1
-q^{k_n}))$ random variable,
and from (\ref{eqNn-expasy}) and (\ref{eqNn-var}) we see that $\mu
_{n,k_n}$ and $\sigma_{n,k_n}^2$ is
the mean\vspace*{1pt} and the variance of the sum, and so (\ref{eqBern}) is exactly
the characteristic function of
a properly centered and normed sum of i.i.d. random
variables. Since $\sigma_{n,k_n} \to\infty$, a simple application of
the Lindeberg--Feller theorem shows that
the limit is $\mathrm{e}^{-t^2/2}$, the characteristic function of the standard
normal distribution. This proves~(\ref{eqNn-norm}).

We turn to the case of typical maximum. For any $j \in\Z$
\[
\bigl( 1 - q^{\lceil\log_r n \rceil+ j - 1} (1-ps) \bigr)^n = \biggl(
1 -
\frac{ \gamma_n q^{j-1} (1-ps)}{n} \biggr)^n \sim \mathrm{e}^{-(1-ps) \gamma_n
q^{j-1} },
\]
and (\ref{eqNn-typ}) follows.

For (\ref{eqNn-large}), it is easy to check that the expectation in
(\ref{eqNn-exp})
tends to 1, whenever $k_n - \log_r n \to\infty$. Since $N_n \geq1$,
the statement follows.
\end{pf*}

For $j \in\Z$ and $m \geq1$, let denote
%
%
\[
\label{eqrj} r_{j,\gamma}(m) = \frac{(p q^{j-1} \gamma)^m}{m!} \bigl
( \mathrm{e}^{pq^{j-1} \gamma}
- 1 \bigr)^{-1}.
\]
Then
\[
h_{j,\gamma} (s) = \sum_{m=1}^\infty
r_{j,\gamma} (m) s^m.
\]
From (\ref{eqNn-typ}), we obtain that
%
%
\begin{equation}
\label{eqr-merge} \lim_{n \to\infty} \max_{1 \leq m \leq n} \bigl
\llvert\p\bigl\{ N_n = m \mid X_n^*= r^{(\lfloor\log_r n\rfloor+ j)/\alpha}
\bigr\} - r_{j, \gamma_n} (m) \bigr\rrvert= 0.
\end{equation}

As a consequence of Proposition \ref{propNn} and Lemma \ref
{lemmamax-conv}, we obtain the unconditional
asymptotic behavior of $N_n$, which also can be described through a
merging phenomenon.

%
\begin{corollary}
Let us denote
\[
h_\gamma(s) = \sum_{j= - \infty}^{\infty}
\bigl( \mathrm{e}^{- ( 1- ps) \gamma q^{j-1} } - \mathrm{e}^{-\gamma q^{j-1}} \bigr).
\]
Then for the generating function of $N_n$ we have
\[
\lim_{n \to\infty} \bigl\llvert g_n(s) -
h_{\gamma_n}(s) \bigr\rrvert= 0, \qquad s \in[0,1].
\]
\end{corollary}

Given that $X \leq r^{k/\alpha}$ for $i \leq k$ we have
$\PROB\{ X = r^{i/\alpha} \mid X \leq r^{k/\alpha} \} =
p q^{i-1} / (1 - q^{k})$. Introduce the
corresponding distribution function
%
%
\begin{equation}
\label{eqcond-df} F_k (x) = \PROB\bigl\{ X \leq x \mid X \leq
r^{k/\alpha} \bigr\} = \cases{ \displaystyle\frac{1}{1- q^{k}} \biggl[
1 -
\frac{ r^{\{\alpha\log_r x\}}}{x^\alpha} \biggr], &\quad for $x \in
\bigl[r^{1/\alpha},
r^{k/\alpha}\bigr]$,
\vspace*{5pt}\cr
1, &\quad for $ x > r^{k/\alpha}$.}
\end{equation}
In the following $X^{(k)}, X_1^{(k)}, \ldots,$ are i.i.d. random
variables with
distribution function $F_k $, and
%
%
\begin{equation}
\label{eqS^k} S_n^{(k)} = X_1^{(k)}
+ \cdots+ X_n^{(k)}
\end{equation}
stands for their partial sums.
By Lemma \ref{lemmacond}, conditioning on $X_n^* = r^{k/\alpha}$, $N_n
= m$
%
%
\begin{equation}
\label{eqcond-repr} S_n \stackrel{\mathcal{D}} {=} m r^{k/\alpha} +
\sum_{i=1}^{n-m} X^{(k-1)}_i
= m r^{k/\alpha} + S_{n-m}^{(k-1)}.
\end{equation}
Calculating the moments we obtain
%
%
\begin{equation}
\label{eqtrunc-moment}
\EXP\bigl(X^{(k)}\bigr)^{\ell} =
\frac{1}{1 - q^{k}} \sum_{i=1}^k
r^{(i \ell)/\alpha} q^{i-1} p 
\cases{ \displaystyle
\frac{p r^{\ell/\alpha}}{1- q^{k}} \frac{r^{ ( {\ell}/{\alpha} - 1
) k} - 1}{r^{{\ell}/{\alpha}-1} - 1},& \quad for $\ell\ne\alpha$,
\vspace*{5pt}\cr
\displaystyle\frac{pr }{1- q^{k}} k, &\quad for $\ell= \alpha$.}
\end{equation}
Note that for $\alpha> \ell$ the truncated $\ell$th moment converges to
$\E X^{\ell}$ as $k \to\infty$, while in other cases the series diverges.

According to Lemma \ref{lemmamax-conv}, the typical values for $X_n^*
= r^{k_n/\alpha}$
are of the form $r^{(\lceil\log_r n \rceil+ j)/\alpha}$,
for some $j \in\Z$.
Therefore, the case $r^{k_n} /n \to0$ corresponds to small maximum,
and $r^{k_n} / n \to\infty$ corresponds to large one. In what follows, we
determine the asymptotic behavior of the sum conditioned on small,
typical and
large maximum.

\subsection{Conditioning on small maximum}

From (\ref{eqcond-repr}), we see that conditioning on the maximum
value $S_n$ is a sum of random number of
i.i.d. random variables. Moreover, (\ref{eqNn-norm}) says that
conditioning on a small maximum
$N_n$ is asymptotically normal. To obtain limit distribution for random
number of i.i.d. random
variables, first we have to determine the behavior of the sum of $n$
i.i.d. random variables.

The following proposition is the conditional counterpart of
Theorem 4 in \cite{GyK} (there only the classical case is treated),
which states that for the sum of truncated variables at $c_n$
the central limit theorem holds if and only if $c_n/n \to0$.
The proof is also similar, therefore we only sketch it.

If we condition on $X_n^* = r^{1/\alpha}$, then all the variables are
degenerate, so we
exclude this case in the following statement. Recall definitions (\ref
{eqcond-df}),
(\ref{eqS^k}) and the notation after it.

%
\begin{proposition} \label{propnormality-uncond}
For $\alpha\in(0,2), k_n \geq2$
%
%
\begin{equation}
\label{eqnormality-uncond} \frac{S_n^{(k_n)} - \EXP S_n^{(k_n)} }{
\sqrt{ \Var ( S_n^{(k_n)} ) } } \stackrel{\mathcal{D}} {\longrightarrow}
N(0,1)
\end{equation}
if and only if $\log_r n - k_n \to\infty$.
\end{proposition}

\begin{pf}
We may assume that $k_n \to\infty$. From equation (\ref
{eqtrunc-moment}), we have
that for any $\alpha\in(0,2)$
%
%
\begin{equation}
\label{eqe-ovar} \bigl( \E X^{(k)} \bigr)^2 = \mathrm{o}\bigl( \E
\bigl(X^{(k)}\bigr)^2 \bigr)\qquad\mbox{as } k \to\infty,
\end{equation}
therefore
\[
\Var X^{(k_n)} \sim\frac{p r^{ 2/\alpha}}{r^{2/\alpha- 1} - 1} r^{ ( 2/\alpha - 1 ) k_n}.
\]
Thus for the variance of the sum
%
%
\begin{equation}
\label{eqvar-sum} s_n^2 = \Var S_{n}^{(k_n)}
= n \Var X^{(k_n)} \sim n \frac{p r^{2/\alpha}}{r^{2/\alpha - 1} - 1} r^{ ( {2}/{\alpha} - 1 ) k_n}.
\end{equation}
By the Lindeberg--Feller central limit theorem
\[
\frac{S_{n}^{(k_n)} - \EXP S_{n}^{(k_n)}}{s_n} \stackrel{\mathcal{D}}
{\longrightarrow} N(0,1)
\]
holds if and only if for every $\varepsilon> 0$
\[
L_n(\varepsilon) = \frac{n}{s_n^2} \int_{\{ \llvert X^{(k_n)} - \EXP
X^{(k_n)} \rrvert > \varepsilon s_n \} }
\bigl( X^{(k_n)} - \EXP X^{(k_n)} \bigr)^2 \,\d\p\to0.
\]
By (\ref{eqe-ovar}), it is easy to show that
\[
L_n(\varepsilon) \sim\frac{n}{s_n^2} \int_{\{ X^{(k_n)} > \varepsilon
s_n \} }
\bigl( X^{(k_n)} \bigr)^2 \,\d\p.
\]
If $r^{k_n} / n \to0$, then by (\ref{eqvar-sum}) the domain of
integration in
$L_n(\varepsilon)$ is empty for large $n$, therefore Lindeberg's
condition holds.

While if $r^{k_n} / n > \varepsilon$ for some $\varepsilon> 0$ and
$n$, then
by (\ref{eqvar-sum}) we have $r^{k_n/\alpha} - \E X^{(k_n)} >
\varepsilon' s_n$ for some $\varepsilon'$, thus
the last jump of $X^{(k_n)}$ belongs to the domain of integration. Therefore,
\[
L_n\bigl(\varepsilon'\bigr) \geq \frac{n}{s_n^2}
r^{{2 k_n}/{\alpha}} q^{k_n - 1} p > \frac{1}{2} \frac{r^{2/\alpha - 1} - 1}{r^{2/\alpha - 1}}.
\]
The proof is complete.\vadjust{\goodbreak}
\end{pf}

Therefore CLT holds for the random index $N_n$ (see (\ref
{eqNn-norm})) and also for the corresponding
deterministic term sums (previous proposition). Combining these two
results the general theory for random sums
(Theorem 4.1.1 in Gnedenko and Korolev \cite{GKor}) implies the following.

%
\begin{proposition} \label{propnormality}
Let $\alpha\in(0,2)$. Given that $X_n^* = r^{k_n / \alpha}$, $k_n
\geq2$, such that $\log_r n - k_n \to\infty$
%
%
\begin{equation}
\label{eqnormality} \frac{S_n - \EXP[ S_n \mid X_n^* = r^{k_n /
\alpha}]}{ \sqrt{ \Var ( S_n \mid X_n^* = r^{k_n / \alpha} ) } }
\stackrel{\mathcal{D}} {\longrightarrow}
N(0,1).
\end{equation}
\end{proposition}

\begin{pf}
By (\ref{eqcond-repr}) given that $X_n^* = r^{k/\alpha}$ we may write
\[
S_n \stackrel{\mathcal{D}} {=} N_n r^{k/ \alpha} +
S_{n- N_n}^{(k-1)} = n r^{k/ \alpha} + \sum
_{i=1}^{n-N_n} \bigl( X_i^{(k-1)} -
r^{k/\alpha} \bigr).
\]
We apply Theorem 4.1.1 in \cite{GKor} to the triangular array
\[
\biggl\{ \frac{X_1^{(k_n - 1)} - r^{k_n / \alpha}}{\sqrt{ \Var
S_n^{(k_n - 1)} }}, \ldots, \frac{X_n^{(k_n - 1)} - r^{k_n / \alpha
}}{\sqrt{ \Var S_n^{(k_n - 1)} }} \biggr
\}_{n \geq1}.
\]
By Proposition \ref{propnormality-uncond}
\[
\sum_{i=1}^n \frac{X_i^{(k_n - 1)} - r^{k_n / \alpha}}{\sqrt{ \Var
S_n^{(k_n - 1)} }} -
\frac{n (\E X^{(k_n -1)} - r^{k_n / \alpha})}{\sqrt{ \Var S_n^{(k_n -
1)} }} \stackrel{\mathcal{D}} {\longrightarrow} N(0,1),
\]
that is condition (1.1) on page~93 in \cite{GKor} holds.
First assume that either $k_n \to k$ for some $k \in\N$, or $k_n \to
\infty$. Put
$u= 1 - \lim_{n\to\infty} \frac{q^{k_n - 1}p }{1- q^{k_n}}$.
Using (\ref{eqNn-var})
%
%
\begin{eqnarray}
\label{eqv} && \lim_{n \to\infty} \frac{r^{k_n / \alpha} - \E X^{(k_n
- 1)} }{\sqrt{ \Var S_n^{(k_n - 1)} }} \sqrt{
\Var\bigl(N_n \mid X_n^* = r^{k_n / \alpha} \bigr)}
\nonumber\\[-8pt]\\[-8pt]\nonumber
&&\quad  = \lim_{n \to\infty} \frac{r^{k_n / \alpha} - \E X^{(k_n - 1)}
}{\sqrt{ \Var X^{(k_n - 1)} }} \sqrt{ \frac{p q^{k_n - 1}}{1- q^{k_n}}
\biggl( 1 - \frac{p q^{k_n - 1}}{1 - q^{k_n}} \biggr) } =: v,
\end{eqnarray}
and the latter limit exists both for $k_n \equiv k$ and for $k_n \to
\infty$.
Using (\ref{eqNn-norm})
\[
\biggl( \frac{n - N_n}{n}, \frac{ n (\E X^{(k_n - 1)} - r^{k_n / \alpha
} )}{\sqrt{ \Var S_n^{(k_n - 1)} }} \frac{n - N_n}{n} -
c_n \biggr) \stackrel{\mathcal{D}} {\longrightarrow} ( u, v Z),
\]
where $Z$ is a standard normal random variable and
\[
c_n = - \bigl( n - \E\bigl[ N_n \mid
X_n^{*} = r^{k_n / \alpha} \bigr] \bigr)
\frac{r^{k_n / \alpha} - \E X^{(k_n -1)} }{\sqrt{ \Var S_n^{(k_n - 1)} }}.
\]
That is, condition (1.9) on page 96 in \cite{GKor} holds,
so Theorem 4.1.1 applies, and we obtain that given $X_n^* = r^{k_n /
\alpha}$
\[
\frac{ \sum_{i=1}^{n-N_n} ( X_i^{(k_n-1)} - r^{k_n/\alpha} ) }{\sqrt{
\Var S_n^{(k_n - 1)} }} - c_n \stackrel{\mathcal{D}} {\longrightarrow}
N\bigl(0,v^2 + u\bigr).
\]

Using (\ref{eqcond-repr}), standard calculation gives that
\[
\E\bigl[ S_n \mid X_n^* = r^{k/\alpha} \bigr] = n \E
X^{(k-1)} + \E\bigl[ N_n \mid X_n^*=
r^{k/\alpha} \bigr] \bigl(r^{k/ \alpha} - \E X^{(k-1)} \bigr),
\]
and
\begin{eqnarray*}
\Var\bigl(S_n \mid X_n^* = r^{k/\alpha} \bigr) &= &
\Var\bigl( N_n \mid X_n^* = r^{k/\alpha} \bigr)
\bigl(r^{k/ \alpha} - \E X^{(k-1)} \bigr)^2
\\
&&{} + \bigl(n - \E\bigl[ N_n \mid X_n^* = r^{k/\alpha}
\bigr] \bigr) \Var X^{(k-1)}.
\end{eqnarray*}
Substituting back the asymptotics (\ref{eqNn-expasy}) and using (\ref
{eqv}) we get that
\[
\lim_{n \to\infty} \frac{\Var S_n^{(k_n -1)} } {\Var(S_n \mid X_n^*
= r^{k_n/ \alpha}) } = \frac{1}{ v^2 + u}.
\]
Summarizing, we obtain (\ref{eqnormality}).

Now let $k_n$ be an arbitrary sequence. From any subsequence $\{ n' \}
$, one can choose a further subsequence
$\{ n'' \}$, such that either $k_{n''} \to k \in\N$ or $k_{n''} \to
\infty$ holds, and so on this subsequence
the convergence takes place. This is equivalent to (\ref{eqnormality}).
\end{pf}

%
\begin{remark}
Without proof, we note that convergence of moments also hold both in
(\ref{eqnormality})
and in (\ref{eqnormality-uncond}).
In view of the distributional convergence, it is enough (in fact equivalent)
to show the uniform integrability of arbitrary powers of the corresponding
random variables.

Using Chernoff's bounding technique, one can prove exponential bounds
for the
tail probabilities
\[
\p\bigl\{ S_n^{(k)} - \E S_n^{(k)} >
n^{1/\alpha} x \bigr\},
\]
from which uniform integrability follows. These bounds and a detailed
proof of the
statement will be published elsewhere, as a continuation of the present paper.
\end{remark}

For $ \alpha> 2$, the variance is finite thus usual central limit
theorem holds without conditioning.
As it was pointed out in the introduction, for $\alpha= 2$ the
generalized St.~$\operatorname{Petersburg}(2,p)$
distribution has infinite variance, but it is still in the domain of
attraction of the normal law.
However, the normalizing sequence is $\sqrt{p r n \log_r n}$, therefore
it is meaningful to ask what is
the necessary and sufficient condition for (\ref{eqnormality-uncond}).
%

%
\begin{proposition} \label{propnormality2}
Let $\alpha= 2$.
Then (\ref{eqnormality-uncond}) holds if and only if
%
%
\begin{equation}
\label{eqnormality-2} \liminf_{n \to\infty} \frac{ \log_r n }{k_n}
\geq1.
\end{equation}
\end{proposition}

Note that the condition is much weaker than the condition for $\alpha
\in(0,2)$. In particular, it also covers the
typical case $k_n \sim\log_r n$, and part of the large maximum case.

\begin{pf*}{Proof of Proposition \ref{propnormality2}}
The proof is exactly the same as in the $\alpha< 2$ case, the only
difference is the
variance asymptotic.

We again assume that $k_n \to\infty$. From equation (\ref
{eqtrunc-moment}), we have for the variance of the sum
%
%
\begin{equation}
\label{eqvar-sum2} s_n^2 = \Var S_{n}^{(k_n)}
= n \Var X^{(k_n)} \sim\frac{p}{q} n k_n.
\end{equation}
By the Lindeberg--Feller theorem, CLT holds if and only if
$L_n(\varepsilon) \to0$ for any
$\varepsilon> 0$. We have
\begin{eqnarray*}
L_n(\varepsilon) & \sim& \frac{n}{s_n^2} \int_{\{ X^{(k_n)} >
\varepsilon s_n \} }
\bigl( X^{(k_n)} \bigr)^2 \,\d\p
\\
& =& \frac{1}{k_n} \bigl\llvert\bigl\{ k\dvt  r^{k/2} > \varepsilon
s_n; k \leq k_n \bigr\} \bigr\rrvert
 = \frac{1}{k_n} \biggl( k_n - \biggl\lfloor
\log_r \frac{\varepsilon^2 p n k_n}{q} \biggr\rfloor\biggr)_+,
\end{eqnarray*}
and the latter goes to 0 if and only if
\[
\liminf_{ n \to\infty} \frac{\log_r (n k_n)}{k_n} \geq1.
\]
Since $(\log_r k_n) / k_n \to0$ this is equivalent to (\ref{eqnormality-2}).
\end{pf*}

\subsection{Conditioning on typical maximum} \label{subsecttyp-max}

According to Lemma \ref{lemmamax-conv}, the typical value for $X_n^*$ is
$r^{(\lceil\log_r n \rceil+ j)/\alpha}$, $j \in\Z$. In the
following, we investigate
this case. Since semistability appears, first we briefly define the
semistable distributions,
and summarize their most important properties. For background, we refer
to Meerschaert and Scheffler \cite{MS} and Megyesi \cite{Megyesi2} and
the references
therein.

Let $Y$ be an infinitely divisible random variable with
characteristic function $\phi(t) = \E( \mathrm{e}^{\i tY})$ in its L\'evy form
(\cite{GK}, page~70), given for each $t\in\R$ by
\[
\phi(t) = \exp\biggl\{ \i t \theta- \frac{\sigma^2 }{2} t^2 + \int
_{-\infty}^0 \beta_t(x) \,\d L(x) + \int
_0^\infty\beta_t(x) \,\d R(x) \biggr
\},
\]
where
\[
\beta_t(x) = \mathrm{e}^{\i t x } - 1 - \frac{\i t x }{{1 + x^2}}.
\]

We describe semistable laws in the present framework as follows:
an infinitely divisible law is semistable if and only if
either it is normal (as a semistable distribution of exponent 2),
or there exist nonnegative bounded functions
$M_L(\cdot)$ on $(-\infty, 0)$ and $M_R(\cdot)$ on $(0, \infty)$, one
of which has strictly positive infimum and the other one either has
strictly positive infimum or is identically zero, such that
$L(x) = M_L(x)/|x|^\alpha$, $x<0$, is left-continuous and
non-decreasing on $(-\infty, 0)$ and $R(x)= -M_R(x)/x^\alpha$, $x>0$,
is right-continuous and non-decreasing on $(0, \infty)$ and
$M_L(c^{1/\alpha}x)= M_L(x)$ for all $x < 0$ and
$M_R(c^{1/\alpha}x) = M_R(x)$ for all $x > 0$, with the same period
$c > 1$.

The following theorem of Kruglov \cite{Kr} highlights the importance of
semi\-stability.
Let $Y_1, Y_2, \ldots$ be independent and identically distributed random
variables with the common distribution function $G$.
If for some centering and norming constants $c_{n_k} \in\R$ and $a_{n_k}>0$
the convergence in distribution
%
%
\begin{equation}
\label{eqsemistable-def} \frac{1}{a_{n_k}} \Biggl( \sum_{j=1}^{n_k}
Y_j -c_{n_k} \Biggr) \stackrel{\mathcal{D}} {
\longrightarrow} W
\end{equation}
holds along a subsequence $\{ n_k \}_{n=1}^\infty\subset\N$
satisfying
%
%
\begin{equation}
\label{eqsemi-c} \lim_{k \rightarrow\infty} \frac{n_{k+1} }{n_k}= c
\qquad
\mbox{for some } c \in[1, \infty),
\end{equation}
then the non-degenerate limit $W$ is necessarily semistable. When the
exponent $\alpha< 2$, the $c$ in the common multiplicative period
of $M_L(\cdot)$ and $M_R(\cdot)$ is the $c$ from the latter growth
condition on
$\{n_k\}$. Conversely, for an arbitrary semistable distribution
there exists a distribution function $G$ for which (\ref{eqsemistable-def})
holds along some $\{ n_k \} \subset\N$ satisfying (\ref{eqsemi-c}).

Now we turn to the asymptotic behavior of $S_n^{(\lfloor\log_r n
\rfloor+ j)}$
defined in (\ref{eqS^k}). Recall the definition of $\mu_n $ in (\ref{eqmu}).

%
\begin{proposition} \label{proptypical-max-uncond}
Let $\alpha\in(0,2)$, $ j \in\Z$. The centered and normed sum
\[
\frac{S_{n_k}^{(\lceil\log_r n \rceil+ j)}}{n_k^{1/\alpha}} - \mu_{n_k}
\]
converges in distribution if and only if $\gamma_{n_k} \to\gamma$, for
some $\gamma\in[q, 1]$.
In this case the limit $W_{j, \gamma} $ has characteristic function
%
%
\begin{equation}
\label{eqchar-f} \varphi_{j,\gamma} (t) = \EXP \mathrm{e}^{\i t W_{j, \gamma} }
= \exp
\biggl[ \i t u_{j, \gamma} + \int_0^\infty
\bigl( \mathrm{e}^{\i t x} - 1 - \i t x \bigr) \,\d L_{j, \gamma} (x) \biggr],
\end{equation}
with
%
%
\begin{equation}
\label{eqL-j} L _{j, \gamma}(x) = \cases{ \displaystyle\gamma
q^j - \frac{r^{ \{ \log_r (\gamma x^\alpha) \} }}{x^\alpha}, &\quad for
$x < r^{j/\alpha}
\gamma^{-1/\alpha}$,
\cr
0, &\quad for $x \geq r^{j/\alpha}
\gamma^{-1/\alpha}$,}
\end{equation}
and
%
%
\begin{equation}
\label{equ} u_{j,\gamma} = \cases{ \displaystyle\frac{p r^{1/\alpha
}}{r^{1/\alpha- 1} - 1}
r^{j (\alpha^{-1} - 1)} \gamma^{1- \alpha^{-1}}, &\quad$\alpha\ne1$,
\vspace*{5pt}\cr
\displaystyle
p r \log_r \frac{r^j}{\gamma}, &\quad$\alpha= 1$.}
\end{equation}
\end{proposition}

Note that the random variables $W_{j,q} $ and $W_{j+1,1} $ have the
same distribution.
This implies that when the set of limit points of the sequence
$\{ \gamma_{n_k} \}_{k \in\N}$ is $\{ q, 1\}$ then convergence in distribution
does not hold, contrary to the unconditional case described after (\ref
{eqmax-merge}).

\begin{pf*}{Proof of Proposition \ref{proptypical-max-uncond}}
Recall the notation in (\ref{eqcond-df}).
According to Theorem 25.1 in Gnedenko and Kolmogorov \cite{GK} the
centered and normalized sum
$S_{n}^{( \lceil\log_r n \rceil+ j)} /n^{1/\alpha} - A_n$ converges in
distribution with some $A_n$ along the subsequence $\{ n_k \}$ if and
only if
%
%
\begin{equation}
\label{eqconv-cond-1} n_k \bigl[1 - F_{\lceil\log_r n_k \rceil+ j}
\bigl(n_k^{1/\alpha}
x\bigr) \bigr] \qquad\mbox{converges}
\end{equation}
and
%
%
\begin{equation}
\label{eqconv-cond-2} n_k F_{\lceil\log_r n_k \rceil+ j} \bigl
(-n_k^{1/\alpha}
x\bigr)\qquad\mbox{converges},
\end{equation}
for any $x > 0$, which is a continuity point of the corresponding limit
function, and
%
%
\begin{eqnarray}
\label{eqconv-cond-3} && \lim_{\varepsilon\to0} \limsup_{k \to\infty}
n_k \int_{|x| \leq\varepsilon} x^2 \,\d
F_{\lceil\log_r n_k \rceil+ j} \bigl( n_k^{1/\alpha} x\bigr)
\nonumber\\[-8pt]\\[-8pt]\nonumber
&&\quad  = \lim_{\varepsilon\to0} \liminf_{k \to\infty}
n_k \int_{|x| \leq\varepsilon} x^2 \,\d
F_{\lceil\log_r n_k \rceil+ j} \bigl( n_k^{1/\alpha} x \bigr) =
\sigma^2.
\end{eqnarray}

Condition (\ref{eqconv-cond-2}) holds for any subsequence with $0$ as
the limit
function.
Using (\ref{eqcond-df}) for $x < r^{j/\alpha}/\gamma_{n_k}^{1/\alpha}$
\begin{eqnarray*}
n_k \bigl[1 - F_{\lceil\log_r n_k \rceil+ j} \bigl(n_k^{1/\alpha}
x\bigr) \bigr] & =& \frac{- n_k q^{\lceil\log_r n_k \rceil+ j} }{ 1 -
q^{\lceil\log_r n_k \rceil+ j}} + \frac{r^{ \{ \log_r (n_k x^\alpha)
\} } x^{-\alpha}}{
1 - q^{\lceil\log_r n_k \rceil+ j}}
\\
& =& \frac{ - q^j \gamma_{n_k} + r^{ \{ \log_r (n_k x^\alpha) \} }
x^{-\alpha} }{
1 - q^{\lceil\log_r n_k \rceil+ j}},
\end{eqnarray*}
thus condition (\ref{eqconv-cond-1}) reduces to the convergence of
\[
- \frac{\gamma_{n_k}}{r^{j}} + \frac{r^{ \{ \log_r ({n_k} x^\alpha) \}
}}{x^{\alpha}}
\]
for $x < r^{j/\alpha}/\gamma_{n_k}^{1/\alpha}$, which is a continuity
point of the limit.
This holds if and only if $\gamma_{n_k}$ converges to some $\gamma\in
[q,1]$, in
which case the limit function is $L_{j, \gamma}$ in
(\ref{eqL-j}), as stated.

Finally, for condition (\ref{eqconv-cond-3}) assume that $\varepsilon
< r^{j/\alpha}$.
Then
\begin{eqnarray*}
n \int_{|x| \leq\varepsilon} x^2 \,\d F_{\lceil\log_r n \rceil+ j}
\bigl(n^{1/\alpha} x\bigr) & =& n^{1-2/\alpha} \int_{|y| \leq
\varepsilon n^{1/\alpha} }
y^2 \,\d F_{\lceil\log_r n \rceil+ j} (y)
\\
& =& n^{1-2/\alpha} \sum_{k\dvt  r^{k/\alpha} \leq\varepsilon
n^{1/\alpha}} r^{2k/\alpha}
\frac{p q^{k-1}}{1- q^{\lceil\log_r n \rceil+ j}}
\\
& \leq&\frac{\varepsilon^{2 - \alpha}}{q - q^{2/\alpha} },
\end{eqnarray*}
for $n$ large enough,
which shows that (\ref{eqconv-cond-3}) holds along the whole sequence
with $\sigma^2 = 0$.

Theorem 25.1 in \cite{GK} states that the centering sequence $A_{n,j}$ can
be chosen as
\[
A_{n,j} = n \int_{|x| \leq\tau} x \,\d F_{\lceil\log_r n \rceil+ j}
\bigl(n^{1/\alpha} x\bigr),
\]
for arbitrary $\tau> 0$. Let us choose $\tau> r^{(j + 1)/\alpha}$. Then
by (\ref{eqtrunc-moment})
\begin{eqnarray*}
A_{n,j} & =& n^{1- \alpha^{-1}} \int_0^{\tau n^{1/\alpha} }
x \,\d F_{\lceil\log_r n \rceil+ j}(x) = n^{1- \alpha^{-1}} \EXP
X^{(\lceil\log_r n \rceil+ j)}
\\
& =& \cases{ \displaystyle\frac{p r^{1/\alpha}}{r^{1/\alpha- 1} -1} r^{
j (\alpha^{-1} - 1)}
\gamma_n^{1-\alpha^{-1}} + \mathrm{o}(1), &\quad$\alpha< 1$,
\vspace*{5pt}\cr
pr \bigl(
\lceil\log_r n \rceil+ j \bigr) + \mathrm{o}(1), &\quad$\alpha= 1$,
\vspace*{5pt}\cr
\displaystyle n^{1- \alpha^{-1}} \E X - \frac{p r^{1/\alpha}}{1 - r^{1/\alpha- 1}}
r^{ j (\alpha^{-1} - 1)}
\gamma_n^{1-\alpha^{-1}} + \mathrm{o}(1), &\quad$\alpha> 1$,}
\end{eqnarray*}
where $\mathrm{o}(1) \to0$ as $n \to\infty$.
We obtain that whenever $\gamma_{n_k} \to\gamma$
\[
\frac{S_{n_k}^{( \lceil\log_r n_k \rceil+ j )}}{n_k^{1/\alpha}} -
A_{n,j} \stackrel{\mathcal{D}} {\longrightarrow}
\widetilde W_{j, \gamma},
\]
where
\[
\EXP \mathrm{e}^{\i t \widetilde W_{j,\gamma} } = \exp\biggl[ \int_0^\infty
\bigl(\mathrm{e}^{\i t x } - 1 - \i t x \bigr) \,\d L_{j, \gamma} (x) \biggr].
\]

Recall the definition of $\mu_n$ in (\ref{eqmu}). We have
\[
\mu_n - A_{n,j} = \cases{ \displaystyle-
\frac{p r^{1/\alpha}}{r^{1/\alpha- 1} - 1} r^{j (\alpha^{-1} - 1)}
\gamma_n^{1- \alpha^{-1}} + \mathrm{o}(1), &
\quad$\alpha< 1$,
\vspace*{5pt}\cr
\displaystyle- pr \bigl( j + \log_r
\gamma_n^{-1} \bigr) + \mathrm{o}(1), &\quad$\alpha= 1$,
\vspace*{5pt}\cr
\displaystyle\frac{p r^{1/\alpha}}{1 - r^{1/\alpha- 1 }} r^{j (\alpha
^{-1} - 1)} \gamma_n^{1- \alpha^{-1}}
+ \mathrm{o}(1), &\quad$\alpha> 1$.}
\]
Therefore,
\[
\frac{S_{n_k}^{( \lceil\log_r n_k \rceil+ j )}}{n_k^{1/\alpha}} - \mu
_{n_k} \stackrel{\mathcal{D}} {
\longrightarrow} \widetilde W_{j, \gamma} + u_{j, \gamma},
\]
with the constant $u_{j,\gamma}$ in (\ref{equ}), as stated.
\end{pf*}

The L\'evy function $L_{j,\gamma} $ is a pure jump function with jumps
at $r^{k/\alpha} \gamma^{-1/\alpha}$, $k \leq j$, such that
$L_{j, \gamma} (r^{k/\alpha} \gamma^{-1/\alpha} ) -
L_{j, \gamma} (r^{k/\alpha} \gamma^{-1/\alpha} - )
= \gamma p q^{k-1}$, for $k \leq j$.
Introduce the notation
\[
G_{j, \gamma} (x)= \PROB\{ W_{j,\gamma} \leq x \}.
\]
The form of the L\'evy function $L_{j, \gamma} $ in
(\ref{eqL-j})
implies that for any $j \in\Z, \gamma\in[q,1]$, the support of
$W_{j,\gamma} $ is $\R$ for $\alpha\geq1$,
while for $\alpha< 1$ the support of $W_{j, \gamma} $ is
$[0, \infty)$, since
\[
u_{j,\gamma} - \int_0^\infty x \,\d
L_{j,\gamma} (x) = 0,
\]
is the drift of the corresponding L\'evy process.
Moreover, the exponential moments $\E \mathrm{e}^{\lambda W_{j, \gamma} }$ are
finite for any $\lambda> 0$, $\alpha\in(0,2)$ and $j \in\Z$,
$\gamma\in[q,1]$,
see, for example, Sato \cite{Sato}, Chapter~5.

The logarithm of the characteristic function of $W_{j, \gamma} $
can be written as
\[
\log\varphi_{j, \gamma} (t) = \i t u_{j, \gamma} + \sum
_{k= - \infty}^j \biggl( \mathrm{e}^{\i t {r^{k/\alpha}}/{\gamma^{1/\alpha
}}} - 1 -
\i t \frac{r^{k/\alpha}}{\gamma^{1/\alpha}} \biggr) \gamma p q^{k-1}.
\]
%
Thus,
\[
\mathfrak{Re} \log\varphi_{j, \gamma} (t) = \sum
_{k=-\infty}^j \biggl( \cos\frac{t r^{k/\alpha}}{\gamma^{1/\alpha}} - 1
\biggr) \gamma p q^{k-1}.
\]
Put
\[
\kappa_{\gamma} (t) = \biggl\lfloor\alpha\log_r
\frac{\gamma^{1/\alpha} \uppi}{2 |t|} \biggr\rfloor.
\]
The same way as in the proof of Lemma 3 in \cite{Cs02} one has that
\begin{eqnarray*}
\mathfrak{Re} \log\varphi_{j, \gamma} (t) & =& - \sum
_{k= - \infty}^j \biggl( 1 - \cos\frac{t r^{k/\alpha}}{\gamma^{1/\alpha}}
\biggr) \gamma p q^{k-1}
\\
& \leq&- \frac{4 p t^2}{q \uppi^2 } \gamma^{1-2/\alpha} \sum
_{k= - \infty}^{ j \wedge\kappa_\gamma(t) } r^{ ( 2/\alpha -
1 ) k}
\\
& \leq&- \frac{4 p \gamma^{1-2/\alpha}}{q \uppi^2 ( 1 - q^{2/\alpha - 1} )} t^2 r^{ ( 2/\alpha - 1 ) (j \wedge\kappa
_\gamma(t) ) }
\\
& \leq&\cases{\displaystyle - c_{\gamma;1} |t|^{\alpha}, &\quad$|t| > T
_\gamma r^{-j/\alpha}$,
\vspace*{3pt}\cr
\displaystyle- c_{\gamma;2} r^{j (2/\alpha - 1 )}
t^2, &\quad$|t| \leq T _\gamma r^{-j/\alpha}$,}
\end{eqnarray*}
where
\begin{eqnarray*}
c _{\gamma;1} &=& \frac{ 2^\alpha p}{
q \uppi^{\alpha} ( r^{2/\alpha - 1} - 1 )},\qquad
c_{\gamma;2} = \frac{4 p \gamma^{1 - 2/\alpha} }{
q \uppi^{2} ( 1 - q^{2/\alpha - 1} )},
\\
T _\gamma &=&  \frac{\gamma^{1/\alpha} \uppi}{2}.
\end{eqnarray*}
By standard Fourier analysis, this implies that $G_{j, \gamma} $ is
infinitely many times differentiable. In particular, by the density
inversion formula
we obtain for $g_{j, \gamma} (x) = (G_{j, \gamma} (x))'$
\begin{eqnarray*}
g_{j,\gamma} ( x) & \leq&\frac{1}{2\uppi} \int
_{-\infty}^\infty\bigl\llvert\varphi_{j, \gamma}
(t)\bigr\rrvert\,\d t
\\
& \leq&\frac{1}{\uppi} \biggl( \int_0^{ T_\gamma q^{j/\alpha}}
\mathrm{e}^{-c_{\gamma;2} r^{j ( 2/\alpha - 1 )} t^2} \,\d t + \int
_{T_\gamma q^{j/\alpha}}^\infty
\mathrm{e}^{- c_{\gamma;1} t^{\alpha} } \,\d t \biggr)
\\
& \leq&\frac{r^{-j (1/\alpha - 1/2 )}}{2 \sqrt{\uppi
c_{\gamma;2} }} + \frac{\Gamma(\alpha^{-1})}{\alpha\uppi(c_{\gamma;1}
)^{1/\alpha}}.
\end{eqnarray*}
Differentiating the characteristic function, we can
compute the first two moments of the variable $W_{j, \gamma} $.
A little calculation gives that
%
%
\begin{equation}
\label{eqW-moments} \EXP W_{j,\gamma} = u_{j, \gamma} \quad\mbox
{and}\quad
\EXP(W_{j,\gamma} )^2 
=
( \EXP W_{j,\gamma} )^2 + \frac{p}{q -q^{2/\alpha}}
\gamma^{1 - 2/\alpha} r^{ ( 2/\alpha - 1 ) j}
\end{equation}
and so
\[
\Var W_{j,\gamma} = \frac{p}{q -q^{2/\alpha}} \gamma^{1 -
2/\alpha}
r^{ ( 2/\alpha - 1 ) j}.
\]
%

As a simple corollary we obtain the following merging theorem.

\begin{corollary} \label{cortypical-max-uncond-merge}
On the whole sequence of natural numbers
%
%
\begin{equation}
\label{eqmax-cond-merge} \sup_{x \in\R} \biggl\llvert\PROB\biggl\{
\frac{S_n^{(\lceil\log_r n \rceil+ j)}}{n^{1/\alpha}} - \mu_n \leq x
\biggr\} - G_{j, \gamma_n} (x)
\biggr\rrvert\to0\qquad\mbox{as } n \to\infty.
\end{equation}
\end{corollary}

\begin{pf}
The simple proof relies upon the same compactness reasoning as the proof
of Theorem 2 in \cite{CsMe02}. We show
that any subsequence $\{ n' \}$ contains a further subsequence on which
(\ref{eqmax-cond-merge}) holds.

Let $\{ n' \}$ be an arbitrary subsequence. The Bolzano--Weierstrass
theorem allows
us to choose a further subsequence $\{ n'' \}$ such that
$\gamma_{n''} \to\gamma$, for some $\gamma\in[q,1]$.
As $\varphi_{j,\gamma_{n''}}(t) \to\varphi_{j,\gamma}(t)$,
by the continuity of $G_{j, \gamma}$ for any $j$ and $\gamma$ we have that
$G_{j, \gamma_{n''}}(x) \to G_{j, \gamma}(x)$ for any $x$.
Using Proposition \ref{proptypical-max-uncond} the statement follows.
\end{pf}

Now we turn to the conditional limit theorem.

%
\begin{proposition} \label{proptypical-max-cond}
For $\alpha\in(0,2)$, $j \in\Z$ we have
%
%
\begin{equation}
\label{eqtypical-max-cond} \sup_{x \in\R} \biggl\llvert\p\biggl\{
\frac{S_n}{n^{1/\alpha} } - \mu_n \leq x \Big| X_n^* =
r^{({\lceil\log_r n \rceil+ j })/{\alpha}} \biggr\} - \widetilde
G_{j, \gamma_n} (x) \biggr\rrvert\to
0,
\end{equation}
where
%
%
\begin{equation}
\label{eqG-tilde} \widetilde G_{j, \gamma} (x) = \sum
_{m=1}^\infty G_{j-1,\gamma} \biggl( x - m
\frac{r^{j/\alpha}}{\gamma^{1/\alpha}} \biggr) r_{j, \gamma} (m).
\end{equation}
\end{proposition}

%
\begin{remark} \label{remarktyp-max-cond}
For any $j \in\Z$ let $(W_{j-1, \gamma} )_{\gamma\in[q,1]}$ be
random variables with characteristic function $\varphi_{j-1,\gamma}$
defined in (\ref{eqchar-f}), and independently let $(M_{j,\gamma
})_{\gamma\in[q,1]}$ be positive integer valued random variables
with generating function $h_{j,\gamma}$
in (\ref{eqh}). Then conditioning on $X_n^* = r^{(\lceil\log_r n \rceil+ j )/\alpha}$ the sum
$\frac{S_n}{n^{1/\alpha} } - \mu_n$ is close in distribution to
%
%
\begin{equation}
\label{eqU{j,gamma}} U_{j,\gamma}= W_{j-1, \gamma} + M_{j, \gamma}
r^{j/\alpha}/\gamma^{1/\alpha}.
\end{equation}
In fact
\[
\p\bigl\{ W_{j-1, \gamma} + M_{j, \gamma} r^{j/\alpha}/
\gamma^{1/\alpha} \leq x \bigr\} = \p\{ U_{j, \gamma} \leq x \} =
\widetilde G_{j,\gamma} (x).
\]
By (\ref{eqh}) $M_{j,\gamma}$ is a Poisson random variable conditioned
on being nonzero, thus
it has finite exponential moments for any $j\in\Z$ and $\gamma\in[q,1]$.
Moreover, $W_{j-1,\gamma}$ and $M_{j,\gamma}$ are independent,
$W_{j-1,\gamma}$ has finite exponential
moments, therefore $U_{j,\gamma}$ also has finite exponential moments.
We can easily determine the moments of $U_{j,\gamma}$. We have
\[
\E M_{j, \gamma} = \sum_{m=1}^\infty m
r_{j,\gamma}(m) = \frac{p q^{j-1} \gamma \mathrm{e}^{p q^{j-1} \gamma}}{\mathrm{e}^{p
q^{j-1} \gamma} -1},
\]
and
\[
\Var M_{j, \gamma} = \E M_{j,\gamma}^2 - ( \E
M_{j,\gamma} ) ^2 = \frac{p q^{j-1} \gamma \mathrm{e}^{p q^{j-1} \gamma}}{\mathrm{e}^{p
q^{j-1} \gamma} -1} - \frac{(p q^{j-1} \gamma)^2 \mathrm{e}^{p q^{j-1} \gamma
}}{(\mathrm{e}^{p q^{j-1} \gamma} -1)^2}.
\]
Therefore, by (\ref{eqU{j,gamma}})
%
%
\begin{eqnarray}
\label{eqU-moments} \E U_{j, \gamma} & =& \E W_{j-1, \gamma} + \E
M_{j, \gamma} \frac{r^{j/\alpha}}{\gamma^{1/\alpha}},
\nonumber\\[-8pt]\\[-8pt]\nonumber
\Var U_{j,\gamma} & =& \Var W_{j-1,\gamma} + \frac{r^{2j/\alpha}}{\gamma
^{2/\alpha}} \Var
M_{j,\gamma}.
\end{eqnarray}
\end{remark}

\begin{pf*}{Proof of Proposition \ref{proptypical-max-cond}}
According to (\ref{eqcond-repr}) conditioning on $X_n^* = r^{(\lceil\log_r n \rceil+ j )/\alpha}, N_n=m$
\[
S_n \stackrel{\mathcal{D}} {=} m r^{(\lceil\log_r n \rceil+ j )/\alpha} +
S_{n-m}^{(\lceil\log_r n \rceil+ j-1)},
\]
and by Corollary \ref{cortypical-max-uncond-merge} we know the
behavior of the latter sum, as for each fixed $m \geq1$
\[
\sup_{x \in\R} \biggl\llvert\p\biggl\{ \frac{S_{n-m}^{(\lceil\log_r
n \rceil+ j-1)}}{n^{1/\alpha} } -
\mu_n \leq x \biggr\} - G_{j,\gamma}(x) \biggr\rrvert\to0.
\]

By the law of total probability
%
\begin{eqnarray}
\label{eqtot-prob} && \p\biggl\{ \frac{S_n}{n^{1/\alpha} } - \mu_n
\leq x \Big|
X_n^* = r^{(\lceil\log_r n \rceil+ j )/\alpha} \biggr\}
\nonumber\\[-8pt]\\[-8pt]\nonumber
&&\quad  = \sum_{m= 1}^n \p\biggl\{
\frac{S_{n-m}^{(\lceil\log_r n \rceil+ j-1)}}{n^{1/\alpha} } - \mu_n
+ \frac{m r^{j/\alpha}}{\gamma_n^{1/\alpha}} \leq x \biggr\} \p
\bigl\{ N_n = m \mid X_n^* = r^{(\lceil\log_r n \rceil+ j )/\alpha} \bigr\}.
\end{eqnarray}
Combining (\ref{eqr-merge}) with (\ref{eqtot-prob}) it is routine to
obtain (\ref{eqtypical-max-cond}).
\end{pf*}

%
\begin{figure}

\includegraphics{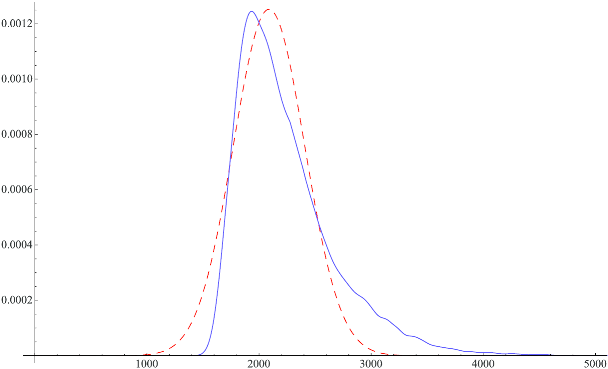}

\caption{The histogram of $S_{n}$ for $n=2^{7}$ ($\alpha= 1$, $p=1/2$)
conditioned on
$X_n^*= 2^{10}$ (solid) and a fitted Gaussian density (dashed).}\vspace*{-6pt}\label{figx2}
\end{figure}

Figure~\ref{figx2} illustrates the histogram of $S_{n}$ for $n=2^{7}$
($\alpha= 1$, $p=1/2$) conditioned on
$X_n^*= 2^{10}$ and a fitted Gaussian density. The histogram has the
property of
positive skewness, which means that the right-hand side tail is larger
than the
left-hand side one. The scaled and translated version of the histogram
corresponds to
the density function of $U_{3,1}$.

\subsection{Conditioning on large maximum}

As we mentioned in the Introduction, the side lobes in Figures~\ref
{fig1} and \ref{fig2}
correspond to the conditional histograms of $\log_2 S_n$ conditioning
the large
values of $X_n^*$, such that
they have disjoint support contained in an interval of length 1. It
means that
$\log_2 X_n^*< \log_2 S_n <\log_2 X_n^* +1$,
or equivalently $X_n^*< S_n <2 X_n^*$, if $X_n^*$ is large enough.
In the next proposition, we make this observation precise.


In the following we investigate the case when $X_n^* = r^{k_n/\alpha}$
is large,
that is,~what happens for $k_n > \log_r n $.
We restrict ourselves to the $\alpha\in(0,2)$ case, since
for $\alpha\geq2$ CLT holds, and thus the corresponding statements
are not interesting.

%
\begin{proposition} \label{proplarge-max}
Let $\alpha\in(0,2)$.
Assume that $k_n - \log_r n \to\infty$. Given that $X_n^* =
r^{k_n/\alpha}$
\[
\frac{S_n}{X_n^{*}} - A_n \stackrel{\p} {\longrightarrow} 1,
\]
where
\[
A_n = \cases{ 0, &\quad$\alpha< 1$,
\vspace*{3pt}\cr
\displaystyle
\frac{p}{q} \frac{n k_n}{r^{k_n}}, &\quad$\alpha= 1$,
\vspace*{3pt}\cr
\displaystyle
\frac{n}{r^{k_n / \alpha}} \frac{p}{q^{1/\alpha} - q}, &\quad$\alpha> 1$.}
\]
\end{proposition}

\begin{pf}
As $k_n - \log_r n \to\infty$ by Proposition \ref{propNn} we have that
$\p\{ N_n = 1 \mid X_n^* = r^{k_n / \alpha} \} \to1$. Therefore, we
may condition
on the event $\{ X_n^* = r^{k_n / \alpha}, N_n = 1 \}$, and given this event
by Lemma \ref{lemmacond}
\[
S_n \stackrel{\mathcal{D}} {=} r^{k_n / \alpha} +
S_{n-1}^{(k_n - 1)}.
\]

Proceeding as in the proof of Proposition \ref
{proptypical-max-uncond}, one can see that
in order to obtain a non-degenerate limit the
normalization for $S_{n-1}^{(k_n -1)}$ should be $n^{1/\alpha}$, but
$r^{k_n/\alpha} / n^{1/\alpha} \to\infty$, so
the maximum alone is too large. That is in this case there is no
non-degenerate limit
distribution.

We shall determine the limit behavior of the
sum $S_{n-1}^{(k_n-1)} / r^{k_n/\alpha} - A_n$, with some
centering $A_n$.
Using Theorem 25.1 in \cite{GK}, one can check as in the proof of
Proposition \ref{proptypical-max-uncond} that condition
(\ref{eqconv-cond-3}) holds, and also (\ref{eqconv-cond-1}) and (\ref
{eqconv-cond-2})
hold with constant 0 as the limit function.
Choosing $\tau> 2$, we get the centering sequence
\[
A_n = n \int_{|x| \leq\tau} x \,\d F_{k_n -1}
\bigl(r^{k_n/\alpha} x\bigr) \sim\frac{n }{r^{k_n/\alpha}} \E X^{(k_n -1)}.
\]

For $\alpha< 1$, using formula (\ref{eqtrunc-moment}) we see that
$A_n \to0$, while for $\alpha= 1$,
\[
A_n \sim\frac{p}{q} \frac{n k_n}{r^{k_n}}.
\]
Finally, for $\alpha> 1$ the expectation $\E X = p/ (q^{1/\alpha} - q)
< \infty$,
therefore
\[
A_n \sim\frac{n}{r^{k_n / \alpha}} \frac{p}{q^{1/\alpha} - q}.
\]

In all cases the limit distribution is degenerate at 0, so we obtain that
%
%
\begin{equation}
\label{eqmax-conv} \frac{S_{n-1}^{(k_n-1)}}{r^{k_n/\alpha}} - A_n \to0,
\end{equation}
in distribution, and so in probability.
Adding the maximum term we obtain the statement.
\end{pf}

%
\begin{remark} \label{remarklarge-max1}
Note that contrary to the case $\alpha\geq1$ for $\alpha< 1$ there
is no need for
centering for any $k_n$ which
satisfies $n/r^{k_n} \to\infty$. That is, given that $X_n^* = r^{k_n/
\alpha}$
\[
\frac{S_n}{X_n^*} \stackrel{\p} {\longrightarrow} 1,
\]
so the maximum term alone dominates the whole sum. This is not
surprising given the
results of Darling \cite{Darling} and Breiman \cite{Breiman}. In
Theorem 5.1
in \cite{Darling} Darling shows that if $Y, Y_1, Y_2, \ldots$ are nonnegative
i.i.d. random variables from the domain of attraction of an $\alpha
$-stable law,
$\alpha\in(0,1)$, then
\[
\frac{\max\{ Y_1, \ldots, Y_n \} }{\sum_{i=1}^n Y_i}
\]
converges in distribution to a non-degenerate random variable. On the
other hand,
Breiman in Theorem 4 \cite{Breiman} proves that this property
characterizes the
domain of attraction. Intuitively, when the tail of the distribution function
behaves as $x^{-\alpha}$, $\alpha\in(0,1)$, the maximum term is about
the same order as the whole sum.
In Proposition \ref{proplarge-max}, we assume that the maximum is
larger than it should
be, so it is reasonable to expect that it dominates the whole sum.

For $\alpha=1$, let us consider the classical case.
For $k_n = \lfloor\log_2 n + \log_2 \log_2 n \rfloor+ j$,
$j \in\Z$, given that $X_n^* = 2^{k_n}$ we again obtain a precise oscillatory
behavior
\[
\frac{S_n}{X_n^*} - 2^{-j} 2^{\{ \log_2 n + \log_2 \log_2 n \}}
\stackrel{\p} {
\longrightarrow} 1.
\]
In fact, (\ref{eqmax-conv}) states more. For
$k_n = \lfloor\log_2 n + a \log_2 \log_2 n \rfloor$, with some $a \in(0,1)$,
given $X_n^* = 2^{k_n}$
\[
\frac{S_n}{X_n^*} - ( \log_2 n )^{1- a} 2^{\{ \log_2 n + a \log_2 \log
_2 n \}}
\stackrel{\p} {\longrightarrow} 1.
\]
Note the interesting phenomenon that although the maximum does not
dominate the sum,
it is large enough to cause a deterministic growth rate.

For $\alpha> 1$ consider the case when $k_n = \lfloor\beta\log_r n
\rfloor$, for
some $\beta> 1$. For $\beta> \alpha$ the centering goes to 0, and so
conditioning
on $X_n^{*} = r^{k_n}$
\[
\frac{S_n}{X_n^*} \stackrel{\p} {\longrightarrow} 1,
\]
thus the maximum dominates the whole sum. For $\alpha= \beta$
we obtain again the oscillatory behavior, as
\[
\frac{S_n}{X_n^*} - \frac{p}{q^{1/\alpha} - q} r^{ {\{ \alpha\log
_r n \}}/{\alpha}} \stackrel{\p} {
\longrightarrow} 1,
\]
while for $1 < \beta< \alpha$ the ratio
grows as $n^{1 - \beta/\alpha} p/(q^{1/\alpha} - q) r^{ \{ \beta\log_r
n \}/\alpha}$.
\end{remark}

%
\begin{remark}
When $A_n = \mathrm{o}(1)$, Proposition \ref{proplarge-max} says that
$S_n / X_n^* \stackrel{\p}{\to} 1$, given $X_n^* = r^{k_n/\alpha}$.
By Chebyshev's inequality one can get the following bound for the rate
of convergence
\[
\PROB\bigl\{ S_n > (1 + \varepsilon) X_n^* \mid
X_n^* = r^{k_n/\alpha} \bigr\} \leq\frac{4pr^{2/\alpha}}{\varepsilon^2
(r^{2/\alpha- 1})}
\frac{n}{r^{k_n}}.
\]
\end{remark}

\section{A series representation of the semistable limit} \label{sectappl}

In this section, $\alpha\in(0,2)$.
The next theorem gives a representation of the semistable
distribution function $G_\gamma$ introduced in (\ref{eqsum-merge}).
Recall the notation $\widetilde G_{j,\gamma}$ in (\ref{eqG-tilde}).
The interesting feature of the statement is that the
distribution functions $\widetilde G_{j,\gamma}$ in the representation
are distribution functions of infinitely divisible random variables with
finite exponential moments. The expectation and variance is
calculated in (\ref{eqU-moments}).

%
\begin{theorem} \label{thmixture}
Let $\alpha\in(0,2)$. For any $\gamma\in[q, 1]$
\[
G_{\gamma} (x) = \sum_{j=-\infty}^\infty
\widetilde G_{j, \gamma} (x) p_{j,\gamma}.
\]
\end{theorem}

%
\begin{remark}
Before the proof we continue Remark \ref{remarktyp-max-cond}. Let
$(W_{j, \gamma})_{j \in\Z, \gamma\in[q,1]}$
be random variables with characteristic function $\varphi_{j,\gamma} $
in (\ref{eqchar-f}),
independently let $(M_{j,\gamma} )_{j \in\Z, \gamma\in[q,1]}$ be
positive integer valued random variables
with generating function $h_{j,\gamma} $ in (\ref{eqh}), and
independently let
$(Y_\gamma)_{\gamma\in[q,1]}$ be integer valued random variable with
probability distribution
$p_{j, \gamma}$. Then
\[
\p\biggl\{ W_{Y_\gamma-1, \gamma} + M_{Y_\gamma, \gamma} \frac
{r^{Y_\gamma/ \alpha}}{\gamma^{1/\alpha}} \leq x
\biggr\} = G_\gamma(x),
\]
or equivalently the semistable random variable $W_\gamma$ has the representation
\[
W_\gamma\stackrel{\mathcal{D}} {=} W_{Y_\gamma-1, \gamma} +
M_{Y_\gamma, \gamma} \frac{r^{Y_\gamma/ \alpha}}{\gamma^{1/\alpha}}.
\]
We note that this probabilistic representation in the classical case
is basically given in Section~8 in~\cite{GML}.
\end{remark}

\begin{pf*}{Proof of Theorem \ref{thmixture}}
We show that for any fixed $x$, one has
\[
\Biggl\llvert\PROB\biggl\{ \frac{S_n}{n^{1/\alpha}} - \mu_n \leq x
\biggr
\}- \sum_{j = -\infty}^\infty\widetilde
G_{j, \gamma_n} (x) p_{j, \gamma_n} \Biggr\rrvert\to0,
\]
which together with formula (\ref{eqsum-merge}) implies the statement.

To ease the notation, introduce
\[
F_{n, j} (x) = \PROB\biggl\{ \frac{S_n}{n^{1/\alpha}} -\mu_n \leq
x \Big| X_n^* = r^{(\lceil\log_r n \rceil+ j)/\alpha} \biggr\}
\]
and
%
%
\begin{equation}
\label{eqq-def} q_{n,j} = \PROB\bigl\{ X_n^* =
r^{(\lceil\log_r n \rceil+ j)/\alpha} \bigr\}.
\end{equation}
By the law of total probability,
\begin{eqnarray*}
&& \PROB\biggl\{ \frac{S_n}{n^{1/\alpha}} - \mu_n \leq x \biggr\}
\\
&&\quad  = \sum_{j= 1 -\lceil\log_r n \rceil}^\infty\PROB\biggl\{
\frac{S_n}{n^{1/\alpha}} - \mu_n \leq x \Big| X_n^* =
r^{(\lceil\log_r n \rceil+ j)/\alpha} \biggr\} \PROB\bigl\{
X_n^* = r^{(\lceil\log_r n \rceil+ j)/\alpha}
\bigr\}
\\
&&\quad  = \sum_{j=1 -\lceil\log_r n \rceil}^\infty F_{n,j} (x)
q_{n,j}.
\end{eqnarray*}
%
For $\varepsilon>0$ choose $j_{\min} < 0 <j_{\max}$
such that for all $n \geq1$
\[
\sum_{j=-\lceil\log_r n
\rceil+1}^{j_{\min}}q_{n,j} <
\varepsilon/4,\qquad\sum_{j=-\infty}^{j_{\min}}p_{j,\gamma_n}
<\varepsilon/4
\]
and
\[
\sum_{j=j_{\max}+1}^{\infty}q_{n,j} <
\varepsilon/4, \qquad\sum_{j=j_{\max}+1}^{\infty}p_{j,\gamma_n}<
\varepsilon/4.
\]
By (\ref{eqmax-merge}) and Lemma \ref{lemmamax-conv}, this is possible.
Thus,
\begin{eqnarray*}
&& \Biggl\llvert\PROB\biggl\{ \frac{S_n}{n^{1/\alpha}} - \mu_n \leq x
\biggr\}- \sum_{j=-\infty}^\infty\widetilde
G_{j,\gamma_n} (x) p_{j,\gamma_n} \Biggr\rrvert
\\
&&\quad \le\sum_{j=-\lceil\log_r n \rceil+1}^{j_{\min}}q_{n,j} +
\sum_{j=-\infty}^{j_{\min}}p_{j,\gamma_n} + \sum
_{j=j_{\max}+1}^{\infty}q_{n,j} + \sum
_{j=j_{\max}+1}^{\infty}p_{j,\gamma_n}
\\
&&\qquad {}  + \Biggl\llvert\sum_{j=j_{\min}+1}^{j_{\max}}
F_{n,j} (x) q_{n,j} - \sum_{j=j_{\min}+1}^{j_{\max}}
\widetilde G_{j, \gamma_n} (x) p_{j,\gamma_n} \Biggr\rrvert
\\
&&\quad \le\varepsilon+ \sum_{j=j_{\min}+1}^{j_{\max}} \bigl
\llvert F_{n,j} (x) - \widetilde G_{j, \gamma_n} (x)\bigr\rrvert+
\sum_{j=j_{\min}+1}^{j_{\max}} \llvert q_{n,j} -
p_{j,\gamma_n} \rrvert\to\varepsilon,
\end{eqnarray*}
where in the last step we applied Lemma \ref{lemmamax-conv} and
Proposition \ref{proptypical-max-cond}.
\end{pf*}

As a consequence of Theorem \ref{thmixture}, using simply Chebyshev's
inequality
combined with the asymptotics of the first and second moments of
$W_{j,\gamma} $
in (\ref{eqW-moments}) one can obtain sharp bounds on the tail of
$G_\gamma$.

%
\begin{corollary}
For any $\gamma\in[q,1]$ for large enough $x$ we have
\[
1 - G_{\gamma} (x) \leq\textnormal{const} \cdot x^{-\alpha}.
\]
\end{corollary}

However, the exact asymptotic behavior of the semistable tail is known.
It follows from a general recent result by Watanabe and Yamamuro \cite{YW}.
Recall that $R_\gamma$ is the L\'evy function of the semistable
limit $W_\gamma$
defined in (\ref{eqLevy-func}). In Theorem 3 in \cite{YW},
they show that
\begin{eqnarray*}
\liminf_{x \to\infty} x^{\alpha} \bigl[1 - G_\gamma
(x)\bigr] & =&  \inf_{1 \leq x \leq r^{1/\alpha}} x^\alpha\bigl(-
R_\gamma(x)\bigr) = 1
\end{eqnarray*}
and
\begin{eqnarray*}
\limsup_{x \to\infty} x^{\alpha} \bigl[1 - G_\gamma
(x)\bigr] & =& \sup_{1 \leq x \leq r^{1/\alpha}} x^{\alpha} \bigl(-
R_\gamma(x-)\bigr) = r.
\end{eqnarray*}

\section*{Acknowledgements}
We are thankful to the anonymous referees for their valuable comments,
which greatly
improved our paper; in particular, for drawing to our attention the paper
by Chow and Teugels \cite{chow-teugels}.

L\'aszl\'o Gy\"orfi was partially supported by the European Union and
the European Social Fund
through project FuturICT.hu (grant no.~T\'AMOP-4.2.2.C-11/1/KONV-2012-0013).
P\'eter Kevei's research
was realized in the frames of T\'AMOP-4.2.4.A/2-11-1-2012-0001
``National Excellence
Program `elaborating and operating an inland student and researcher
personal support
system'. The project was subsidized by the European Union and
co-financed by the
European Social Fund.''
P\'eter Kevei's research was partially supported by the Hungarian Scientific
Research Fund
OTKA PD106181. Dedicated to the memory of S\'andor Cs\"org\H{o}.




\printhistory
\end{document}